# Symbolic Computation of the Roots of any Polynomial with Integer Coefficients

**Ashok Kumar Mittal**
**Department of Physics,**
**Allahabad University, Allahabad – 211 002, India**
**(Email address: mittal_a@vsnl.com)**

**Ashok Kumar Gupta**
**Department of Electronics and Communication,**
**Allahabad University, Allahabad - 211 002, India**
**(Email address: akgjkiapt@hotmail.com)**

**Abstract:** The roots of any polynomial of degree $m$ with integer coefficients, can be computed by manipulation of sequences made from $2m$ distinct symbols and counting the different symbols in the sequences. This method requires only 'primitive' operations like replacement of sequences and counting of symbols. No calculations using 'advanced' operations like multiplication, division, logarithms etc. are needed. The method can be implemented as a geometric construction of roots of polynomials to arbitrary accuracy using only a straight edge, a compass, and pencils of $2m$ different colors. In particular, the ancient problem of the "doubling of cube" is soluble asymptotically by the above-mentioned construction. This method, by which a cube can be doubled, albeit, in infinite steps, is probably the closest to the original problem of construction using only a straight edge and compass in a finite number of steps.

Moreover, to every polynomial of degree $m$ over the field of rationals, can be associated an $m$-term recurrence relation for generating integer sequences. A set of $m$ such sequences, which together exhibit interesting properties related to the roots of the polynomial, can be obtained if the $m$ initial terms of each of these $m$ sequences is chosen in a special way using a matrix associated with the polynomial. Only two of these integer sequences need to be computed to obtain the real root having the largest absolute value. Since this method involves only integers, it is faster than the conventional methods using floating-point arithmetic.

## 1. Introduction

Rational approximates to the largest (absolute value) real root of any polynomial of degree $m$ with integer coefficients can be obtained to arbitrary accuracy by a method requiring only operations consisting of replacement of sequences and counting of symbols. The replacement rules for this method are obtained from a matrix, called the 'replacement matrix', whose elements are related to the coefficients of the polynomial.

The eigenvalues of the 'replacement matrix' are the roots of the polynomial multiplied by the coefficient of the highest power of the polynomial. The largest eigenvalue can be obtained by the method of matrix iteration, provided the largest eigenvalue is unique [1]. It is possible to carry out this matrix iteration method using only 'primitive operations' like replacement of sequences and counting of symbols.

The method can be extended, by modifying the 'replacement matrix', to obtain rational approximates to other roots of the polynomial.

The method can be implemented geometrically using only a straightedge, compass, and pencils of $2m$ different colors. In particular, by this geometrical implementation, it is possible to obtain a length segment whose cube will have a volume arbitrarily close to twice the volume of a given cube, using only a straight edge and compass. A classical problem of antiquity, the doubling of a cube [6], also known as the Delian problem, required construction by straight edge and compass only in finite number of steps. It is known that this is impossible. By bending the rules of permitted construction, for example under Neusis construction, which permits sliding of a marked ruler, it is possible to solve the problem. The method of this paper, by which a cube can be doubled, albeit, in infinite steps, is perhaps the closest to the original problem.

Further, to any polynomial of degree $m$ over the field of rational numbers, one can associate an $m$-term recurrence relation for generating integer sequences. A set of $m$ such sequences, which together exhibit interesting properties related to the roots of the polynomial, can be obtained if the initial $m$ terms of each of these $m$ sequences are chosen in a special way using a matrix related to the polynomial.

## 2. Replacement Matrix and the Replacement Rules

Let $p(x) = a_0 x^m + a_1 x^{m-1} + \ldots + a_m$, where $a_i$ are integers, be the given polynomial with integer coefficients. The 'replacement matrix' corresponding to the polynomial is defined by

$$\mathbf{R} = \begin{bmatrix} -a_1 & -a_2 & -a_3 & \ldots & -a_{m-1} & -a_m \\ a_0 & 0 & 0 & \ldots & 0 & 0 \\ 0 & a_0 & 0 & \ldots & 0 & 0 \\ . & . & . & \ldots & . & . \\ . & . & . & \ldots & 0 & . \\ 0 & 0 & 0 & \ldots & a_0 & 0 \end{bmatrix} \quad (1)$$

Let $r_i$, $i = 1, 2, \ldots, m$ be the roots of $p(x)$. Then the eigenvalues of $\mathbf{R}$ are $a_0 r_i$ with corresponding eigenvector $[r_i^{m-1}, r_i^{m-2}, \ldots, r_i, 1]^T$.

Let $A = \{A_1, A_2, \ldots, A_m\}$ be a finite set of $m$ symbols. Elements of A may be called 'letters' belonging to the 'alphabet' A. Let $A^p = A \times A \times \ldots \times A$ (p times). Elements of $A^p$ are called 'words' of length p that can be formed from 'letters' of

A. Let $A^* = \bigcup_p A^p$. Clearly, $A^*$ consists of all 'words' that can be made from the 'alphabet' A. The replacement rule, that is used to manipulate the symbolic sequences consisting of these symbols, is obtained from the 'replacement matrix' **R**. Let us first consider the case when all the elements of this matrix are non-negative. Let $R_{ij}$ denotes the element in the $i^{th}$ row and $j^{th}$ column of the matrix **R**. The replacement rule R, gives a rule to replace a word W by a new word $W^*$. This rule consists of replacing each letter $A_j$ in W by a word consisting of the $i^{th}$ symbol repeated $R_{ij}$ times. Thus,

$$R(A_j) = (A_1)\wedge R_{1j}\ (A_2)\wedge R_{2j} \ldots\ldots (A_i)\wedge R_{ij} \ldots\ldots (A_m)\wedge R_{mj} \qquad (2)$$

where $(A_i)\wedge k = A_i^k$ stands for repeating k times the symbol $A_i$. This replacement rule fails to have any meaning if any element of the matrix **R** is negative.

In order to find a replacement rule applicable for the case when some of the elements of the 'replacement matrix' are negative integers, we enlarge the 'alphabet' A to have $2m$ symbols. We denote this enlarged 'alphabet by $A = \{A_1, A_2, \ldots, A_m, A_{m+1}, \ldots, A_{2m}\}$. The symbols $A_{m+i}$ are also denoted by $\tilde{A}_i$. The symbols $A_i$ and $\tilde{A}_i$ will be called conjugate to each other so that $A_i^c = \tilde{A}_i$ and $(\tilde{A}_i)^c = A_i$. Then $(A_i)^{-k}$ stands for repeating k times the symbol $A^c$. The replacement rule (2) will have to be extended to provide for replacement of $2m$ symbols. This is given by

$$R(A_j) = (A_1)\wedge R'_{1j}\ (A_2)\wedge R'_{2j} \ldots\ldots (A_i)\wedge R'_{ij} \ldots\ldots (A_{2m})\wedge R'_{(2m)j} \qquad (3)$$

where $R_{ij}$ are now matrix elements of the $2m \times 2m$ extended matrix

$$\mathbf{R_e} = \begin{bmatrix} R & 0 \\ 0 & R \end{bmatrix} \qquad (4)$$

Let $W = A_{s1}A_{s2}\ldots..A_{sq}$, where $A_{sk}$ (k=1,2,........,q) $\in$ A. Then $W \in A^*$. The replacement rule R induces a mapping $R^*: A^* \to A^*$ defined by

$$\begin{aligned} W^* &= R^*(W) \\ &= R^*(A_{s1}A_{s2}\ldots..A_{sq}) \\ &= R(A_{s1})R(A_{s2})\ldots..R(A_{sq}) \end{aligned} \qquad (5)$$

Let $W_0 \in A^*$. Then

$$W_i = (R^*)^i (W_0) = R^*((R^*)^{i-1}(W_0)) = R^*(W_{i-1}) \qquad (6)$$

denotes the word obtained by i times repeated application of the replacement rule (3) on the initial word $W_0$.

Let $\mathbf{n}: A^* \to V^m$ be a mapping which assigns to a word in $A^*$ a vector

$$\mathbf{n}(W) = (n_1(W) - n_1\tilde{\ }(W), \ldots, n_i(W) - n\tilde{\ }_i(W), \ldots, n_m(W) - n\tilde{\ }_m(W))^T \qquad (7)$$

in an m-dimensional vector space $V^m$ such that $n_i(W)$ and $n\tilde{\ }_i(W)$ are non-negative integers denoting the number of times that $A^+_j$ and $A\tilde{\ }_i$ occur, respectively, in the word W. It is straightforward to see that

$$\mathbf{n}(W^*) = \mathbf{n}(R^*(W)) = \mathbf{R}(\mathbf{n}(W)) \qquad (8)$$

Equations (6) and (8) imply

$$\mathbf{n}(W_i) = \mathbf{n}(R^*(W_{i-1})) = \mathbf{R}(\mathbf{n}(W_{i-1})) = \mathbf{R}^i(\mathbf{n}(W_0)) \qquad (9)$$

For almost any $W_0$, $\mathbf{n}(W_i)$ tends, as $i \to \infty$, to multiples of the eigenvector of $\mathbf{R}$ corresponding to the eigenvalue with maximum absolute value [1], provided it is unique. Thus one finds that for almost any initial word $W_0$,

$$\lim_{i \to \infty} n_j(W_i)/n_{j+1}(W_i) = \lambda \qquad j = 1, 2, \ldots, m-1 \qquad (10)$$

The operation of subtraction in equation (7) can be carried out by removing pairs of symbols $A_i$ and $A\tilde{\ }_i$ from each replacement sequence. Thus, Eqn (10) shows that the largest real root of any polynomial with integer coefficients can be obtained to any accuracy only by operations consisting of replacement of sequences and counting of symbols. This method does not require any 'advanced' operation such as multiplication, division, logarithms and anti-logarithms, although the justification of the method uses 'advanced' concepts like the theory of matrices.

### 3. Example

Let

$$p(x) = 2x^2 + 2x - 1 \qquad (11)$$

be the given polynomial. Here $m = 2$, $a_0 = 2$, $a_1 = 2$ and $a_2 = -1$. The 'replacement matrix' corresponding to this polynomial is given by

$$\mathbf{R} = \begin{bmatrix} -2 & 1 \\ 2 & 0 \end{bmatrix} \qquad (12)$$

The extended 'replacement matrix' is given by

$$\mathbf{R}_e = \begin{bmatrix} -2 & 1 & 0 & 0 \\ 2 & 0 & 0 & 0 \\ 0 & 0 & -2 & 1 \\ 0 & 0 & 2 & 0 \end{bmatrix} \quad (13)$$

Using the alphabet A = {0,1,0˜,1˜}, this 'replacement matrix' gives rise to the following replacement rules:

$$\begin{aligned} 0 &\rightarrow 0˜0˜11 \\ 1 &\rightarrow 0 \\ 0˜ &\rightarrow 001˜1˜ \\ 1˜ &\rightarrow 0˜ \end{aligned} \quad (14)$$

This replacement rule generates the following sequence of sequences:

0

0˜0˜11

001˜1˜ 001˜1˜ 0 0

0˜0˜11 0˜0˜11 0˜ 0˜ 0˜0˜11 0˜0˜11 0˜ 0˜ 0˜0˜11 0˜0˜11

001˜1˜ 001˜1˜ 0 0 001˜1˜ 001˜1˜ 0 0 001˜1˜ 001˜1˜ 001˜1˜ 001˜1˜ 0 0 001˜1˜ 001˜1˜ 0 0 001˜1˜ 001˜1˜ 001˜1˜ 001˜1˜ 0 0 001˜1˜ 001˜1˜ 0 0

The number count of different symbols in the above sequences can be summarized in the following table:

| n(0) | n(1) | n(0)/n(1) |
|---|---|---|
| 1 | 0 | inf |
| -2 | 2 | -1 |
| 6 | -4 | -1.5 |
| -16 | 12 | -1.25 |
| 44 | -32 | -1.375 |
| -120 | 88 | -1.36364 |

| 328 | -240 | -1.36667 |
| -896 | 656 | -1.36585 |
| 2448 | -1792 | -1.36607 |

As expected the ratio **n**(0)/**n**(1) converges to the largest (absolute value) root, namely $(-1 - \sqrt{3})/2$, of polynomial p(x) in (11).

## 4. Obtaining some other roots

The eigenvalues of the matrix $\mathbf{R'} = (\alpha \mathbf{I} + \beta \mathbf{R})$ are $r_i' = (\alpha + \beta a_0 r_i)$ with corresponding eigenvectors $[r_i^{m-1}, r_i^{m-2}, \ldots, r_i, 1]^T$, where **I** is the *m x m* identity matrix. A different set of replacement rules may be obtained by replacing matrix **R** in (1) by **R'**, where $\alpha, \beta$ are suitable integers, chosen so that some other root of p goes to the largest (absolute value) eigenvalue of the **R'** without change in its corresponding eigenvector. In the above example, the roots of p are $(-1 - \sqrt{3})/2$ and $(-1 + \sqrt{3})/2$. The largest (absolute value) root is $(-1 - \sqrt{3})/2$. If we take $\alpha = 2$ and $\beta = 1$, we get **R'** given by

$$\mathbf{R'} = \begin{bmatrix} 0 & 1 \\ 2 & 2 \end{bmatrix} \quad (15)$$

The eigenvalues of **R'** are $(1 - \sqrt{3})$ and $(1 + \sqrt{3})$. The second eigenvalue now becomes the largest (absolute value) but the eigenvector remains unchanged.

The replacement rule corresponding to this matrix is given by

0 → 1 1
1 → 0 1 1

This replacement rule generates the following sequence of sequences:

|  | **n**(0) | **n**(1) | **n**(0)/**n**(1) |
|---|---|---|---|
| 0 | 1 | 0 | inf |
| 11 | 0 | 2 | 0 |
| 011 011 | 2 | 4 | 0.5 |
| 11 011 011 11 011 011 | 4 | 12 | 0.3333 |
| 011 011 11 011 011 11 011 011 011 011 11 011 011 11 011 011 | 12 | 32 | 0.375 |

The ratio **n**(0)/**n**(1) converges to the other root, namely ( -1 + √3)/2, of polynomial p(x) in (9). In this way, the method can be modified to obtain some other roots.

## 5. Geometrical implementation of the method

One can implement the above method as a geometrical construction of a line segment, the ratio of whose length to that of a unit length marked on a straight edge, is equal to a root of the polynomial. This construction uses only the straight edge, a compass, and pencils of at most 2*m* different colors. This can be illustrated by considering the famous "doubling the cube" problem. For this one requires pencils of just three colors, say blue (B), green (G) and red (R). Let the side of the cube, which is to be doubled, be called the unit length. One begins by using the straight edge to draw a line segment of unit length with the blue pencil. The method of section 3, modified by section 4, applied to the polynomial p(x) = $x^3$ - 2 leads to the replacement rule [3,5]

B → BG
G → GR
R → RBB

Using this replacement rule, below the first line, one draws a unit line segment with blue pencil, followed by a unit line segment with green pencil. Following the replacement rule, below this line one draws line segments in the sequence BGGR. In the next line one draws BGGRGRRBB. This process is continued an arbitrary number of times. Thereafter, on the next line one first draws a unit green segment for every unit green segment in the last line one after another. This is followed by red segments in the same way, so that on the last line we have a green line followed by a red line.

In the line so drawn, we mark the beginning of the green line as the point A, the junction point between green and red as the point B and the end of the red line as the point C. Inclined to the line AC, one draws a line CD of unit length. The point D is joined to B. With the help of ruler and compass, a line is drawn from A, parallel to BD. Let this line intersect CD produced at the point E. Now if a cube is built with side DE, its volume will be approximately equal to twice the volume of the original cube. Arbitrary accuracy can be obtained by increasing the number of times that the replacement rule is applied.

## 6. Recurrence Relations

The matrix $a_0^{-1}\mathbf{R}$ satisfies the polynomial p(x). Hence,

$$a_0 (a_0^{-1}\mathbf{R})^m + a_1(a_0^{-1}\mathbf{R})^{m-1} + \ldots\ldots + a_m \mathbf{I} = 0 \qquad (16)$$

Multiplying by $a_0^{m-1}$,

$$R^m + a_1 R^{m-1} + a_0 a_2 R^{m-2} + \ldots\ldots + a_0^{m-1} a_m I = 0 \tag{17}$$

Let

$$S_0 = (S^{(1)}_0, S^{(2)}_0, \ldots\ldots, S^{(m)}_0)^T \tag{18}$$

where $S^{(i)}_0$ are arbitrary integers. Define

$$S_j = R^j S_0, \quad j = 1, 2, \ldots\ldots, m-1 \tag{19}$$

The first $m$ terms of the sequence $S^{(i)}$ are taken to be $\{S^{(i)}_0, S^{(i)}_1, \ldots\ldots, S^{(i)}_{m-1}\}$. Beyond this the terms of the sequence are given by the $m$-term recurrence relation

$$S^{(i)}_j = - a_1 S^{(i)}_{j-1} - a_0 a_2 S^{(i)}_{j-2} - \ldots - a_0^{m-1} a_m S^{(i)}_{j-m}, \quad j > m - 1, i = 1, 2, \ldots, m \tag{20}$$

These sequences satisfy the interesting property that each of the ratios $S^{(i)}_j$ to $S^{(i+1)}_j$, $i = 1, 2, \ldots\ldots, m-1$ tends to the root of p, having the largest absolute value, if it is unique.

For our example of equation (11), $m = 2$, $a_0 = 2$, $a_1 = 2$ and $a_2 = -1$.

We take

$$S_0 = [1 \ 0]^T \tag{21}$$

Then

$$S_1 = R S_0 = [-2 \ 2]^T \tag{22}$$

Hence the first two terms of the sequence $S^{(1)}$ are $\{1, -2\}$ and the first two terms of the sequence $S^{(2)}$ are $\{0, 2\}$. Subsequent terms in these sequences are to be obtained by the recursion relation

$$S^{(i)}_j = -2 S^{(i)}_{j-1} + 2 S^{(i)}_{j-2}, \quad j > 1, i = 1, 2 \tag{23}$$

The application of this recursion leads to the following sequences

| j | $S^{(1)}_j$ | $S^{(2)}_j$ | Ratio |
|---|---|---|---|
| 0 | 1 | 0 | inf |
| 1 | -2 | 2 | -1 |
| 2 | 6 | -4 | -1.5 |
| 3 | -16 | 12 | -1.25 |
| 4 | 44 | -32 | -1.375 |
| 5 | -120 | 88 | -1.3636 |
| 6 | 328 | -240 | -1.3667 |

The ratio $S^{(1)}_j$ to $S^{(2)}_j$ converges to the largest (absolute value) root, namely $(-1 - \sqrt{3})/2$, of the polynomial p(x) in (11).

To obtain the other root, one uses the matrix **R'** of (15) in eqn (22) to obtain $S_1$ from an arbitrary $S_0$. Comparing with (17), it is noted that the matrix **R'** satisfies the equation $(R' - 2 I)^2 + 2 (R' - 2 I) - 2 I = 0$ to yield

$$R'^2 - 2 R' - 2 I = 0 \qquad (24)$$

This gives rise to the recurrence relation

$$S^{(i)}_j = 2 S^{(i)}_{j-1} + 2 S^{(i)}_{j-2}, \qquad j > 1, i = 1, 2 \qquad (25)$$

This recurrence relation leads to two sequences, the ratio of whose corresponding terms converges to the other root of the polynomial.

## 7. Conclusion

We have adapted the method of matrix iteration for finding the roots of polynomials to obtain two methods for finding rational approximates to real roots of polynomials. The first method, though impractical, is interesting because it requires only operations of replacement of sequences and counting of symbols. It does not need any 'advanced operation' like multiplication, division, logarithms and anti-logarithms, although the justification of the method uses 'advanced' concepts like the theory of matrices. Numbers that can be obtained using only a ruler and a compass are called constructible numbers [2]. Every constructible number belongs to some $Q(\sqrt{m_1}, \sqrt{m_2}, ....., \sqrt{m_n})$, an extension field of rational numbers by square roots. Conversely, every number belonging to $Q(\sqrt{m_1}, \sqrt{m_2}, ....., \sqrt{m_n})$ is a constuctible number. A number that may be obtained as a root of a polynomial by our replacement of symbols method may be called *constructible by symbolic computation.* It will be interesting to characterize the set of numbers that are *constructible by symbolic computation.* Our method may be implemented geometrically using only a straight edge, compass, and pencils of 2*m* different colors. In particular, doubling of cubes is possible by the above-mentioned construction. This method, by which a cube can be doubled, albeit, in infinite steps, is probably the closest to the original problem of doubling the cube using only straight edge and compass in finite number of steps.

The second method is based on integer arithmetic only. Further, it needs to compute two sequences only with the help of a recurrence relation. This makes the method faster than the conventional methods based on floating point arithmetic or the unmodified matrix iteration method.

This paper extends our earlier work related to integer monic polynomials [3,4,5]. The method can be extended to obtain complex roots. The details will be given elsewhere.

## 8.  Bibiliography